\documentclass[a4paper]{amsart}
\usepackage{algorithm}
\usepackage{algorithmic}
\usepackage{amsfonts}

\author{Konstantinos A. Blekos}
\title{Unranking permutations in transposition order and linear time.}
\date{August 2007}

\begin{document}
\begin{abstract}
	An algorithm is presented for unranking permutations in transposition order:
	Given a seed $s\in\mathbb{N}$, the algorithm produces a permutation $\mathcal{P}(s)$ that
	differs  from the permutation $\mathcal{P}(s+1)$ by the transposition of two
	elements.
\end{abstract}
\maketitle
\tableofcontents

\section{Introduction}
An algorithm $\mathcal{P}(s)$ is constructed that produces a permutation of
objects in factoradic representation. All permutations of $n$ objects can be
produced, evaluating $\mathcal{P}(s)$ for every $0\le s\le n!-1$.

The key property of the algorithm is that permutations $\mathcal{P}(0),
\mathcal{P}(1), \ldots$ are produced in transposition order.

\subsection{Transposition Order}
It is an ordering of permutations in which each two adjacent
permutations differ by the transposition of two elements. For the permutations
of $\{1,2,3\}$ there are two listings which are in transposition order. One is 123,
132, 312, 321, 231, 213, and the other is 123, 321, 312, 213, 231, 132.\cite{wolf}

\subsection{Factoradic representation}
Factoradic is a numeral system based on factorials. In factoradic, the sequence $a_n a_{n-1} \cdots a_0$ represents the number \[\sum_{i=0}^n a_i i!\], where $a_i \le i$

One can represent any permutation using factoradic as follows: Suppose we have
a set of $n+1$ objects $\{k_0, \ldots, k_n\}$ and the factoradic sequence $a_n
a_{n-1} \cdots a_0$.  We remove the $a_n$-th element from the set of objects
(the element with index $a_n$) and place it first on the permuted list. We
continue with the $a_{n-1}$-th object of the re-indexed set, placing it second
on the permuted list and proceed likewise through $a_0$-th (which is always
zero: the first and only element left in the set).\cite{wik}

\section{Algorithm}
The algorithm is $O(n)$, i.e. linear in respect to the number of objects $n$ in the set.

For a given seed $s$, we choose an $n$ such that $s\le n!-1$. The algorithm
outputs the factoradic sequence $f_{n-1} f_{n-2} \cdots f_0$ that corresponds to  the 
$s$-permutation of a set of $n$ objects.

\begin{algorithm}
	\caption{Given seed $s\in\mathbb{N}$ calculate permutation $f_{n-1}f_{n-2}\ldots f_1$}
	\begin{algorithmic}[1]
		\REQUIRE $n$ such that $s\le n!-1$
		\STATE $d_{n+1} = 0$
		\FOR{$k=n$ to $1$}
		\STATE $x_k = \lfloor \frac{s\pmod{k!}}{(k-1)!}\rfloor$
		\STATE $d_k = \biggl\lfloor \frac{\lfloor\frac{s+d_{k+1}(k+1)!}{k!}\rfloor\pmod{(k+1)^2}}{k+2} \biggr\rfloor$
		\STATE $f_{k-1} = (x_k-\lfloor\frac{s}{k!}\rfloor-d_k)\pmod{k}$
	\ENDFOR
\end{algorithmic}
\end{algorithm}

\subsection{Properties}
Defining a distance $d(s,s')$ as the minimum number of transpositions needed to reach
permutation $\mathcal{P}(s')$ starting from permutation
$\mathcal{P}(s)$\footnote{$d(s,s')$ is a proper distance function since
\begin{enumerate}
	\item $d(s,s')\geq 0$
	\item $\ d(s,s')=0\Leftrightarrow \mathcal{P}(s)=\mathcal{P}(s')$
	\item $d(s,s') =d(s',s)$ 
	\item $d(a,b) \leq d(a,c)+d(c,b)$ 
\end{enumerate}
}
, the
following properties can be proved:


\begin{itemize}
	\item $\mathcal{P}$ is a bijection: $ \mathcal{P}(s) \rightarrow f_{n-1}f_{n-2}\ldots f_0 \Leftrightarrow f_{n-1}f_{n-2}\ldots f_0\rightarrow \mathcal{P}(s), \forall s\in\mathbb{N}, \forall f_{n-1}f_{n-2}\ldots f_0\in\mathbb{P}^n $ 
	\item $\mathcal{P}$ unranks permutations in transposition order: $d(s\pm 1, s) = 1$ 
\end{itemize}

Also, 
\begin{itemize}
	\item $d(s,s') \le \min(|s-s'|,n-1)$
		$\Rightarrow d(s\pm 2, s) = 2$
	\item $d(s,0) \leq k-1,\ s < k!$
	\item $|d(s,0)-d(s',0)|\leq d(s,s')$
\end{itemize}

\section{Examples}
For $s = 4$ we choose $n=4$, so:
\begin{eqnarray*}
	x_4 &= &\lfloor\frac{4\pmod{4!}}{3!}\rfloor = 0\\
	d_4 &= &\biggl\lfloor \frac{(\frac{4+(0)(5!)}{4!})\pmod{5^2}}{6} \biggr\rfloor = 0\\
	f_3 &= &(0-\lfloor\frac{4}{4!}\rfloor-0)\pmod{4} =0\\
	x_3 &= &\lfloor\frac{4\pmod{3!}}{2!}\rfloor = 2\\
	d_3 &= &\biggl\lfloor \frac{(\frac{4+(0)(4!)}{3!})\pmod{4^2}}{5} \biggr\rfloor = 0\\
	f_2 &= &(2-\lfloor\frac{4}{3!}\rfloor-0)\pmod{3} =2\\
	\ldots\\
	f_1 &= &0\\
	f_0 &= &0\\
\end{eqnarray*}
thus $\mathcal{P}(4) \rightarrow 0_3 2_2 0_1 0_0\rightarrow 1 4 2 3$.

Similarly:
\begin{eqnarray*}
 \mathcal{P}(5) \rightarrow 0_3 2_2 1_1 0_0 \rightarrow 1 4 3 2\\
 \mathcal{P}(6) \rightarrow 1_3 2_2 1_1 0_0 \rightarrow 2 4 3 1\\
 \mathcal{P}(7) \rightarrow 1_3 2_2 0_1 0_0 \rightarrow 2 4 1 3\\
 \mathcal{P}(319) \rightarrow  2_5 1_4 2_3 2_2 0_1 0_0 \rightarrow 325614\\
 \mathcal{P}(320) \rightarrow  2_5 1_4 2_3 0_2 1_1 0_0 \rightarrow 325164\\
 \mathcal{P}(321) \rightarrow  2_5 1_4 2_3 0_2 0_1 0_0 \rightarrow 325146\\
 \mathcal{P}(322) \rightarrow  2_5 1_4 2_3 1_2 0_1 0_0 \rightarrow 325416\\
\end{eqnarray*}

Using the above values, we can validate some distance relationships: 
\[d(4,5) = 1\]\[d(4,6) = 2\]
\begin{eqnarray*}
	\left|d(319,0) - d(5,0)\right| & = & |3-1| = 2 \\
	& \leq & d(319,5)  = 2 \\
	& \leq & d(319,0) + d(0,5) = 3 + 1 = 4
\end{eqnarray*}
\begin{eqnarray*}
	d(319,5) & =  &2 \\
	& \leq & \min(|319 - 5| , 5) = 5
\end{eqnarray*}


\end{document}